\documentclass[reqno]{amsart}

\makeatletter
\long\def\@makefntext#1{%
  \noindent\@makefnmark#1%
}
\makeatother

\usepackage{amsmath, amssymb, amsthm}
\usepackage{mathtools}
\usepackage{mathrsfs}
\usepackage{bbm}
\usepackage{cases}

\usepackage{graphicx} 
\usepackage{svg}
\usepackage{booktabs}
\usepackage{siunitx}
\usepackage{float} 
\usepackage{subcaption}

\usepackage{url}
\usepackage{enumitem}
\usepackage{color}
\usepackage{epstopdf}

\usepackage{hyperref} 

\usepackage{amsrefs}

\newtheorem{theorem}{Theorem}[section]
\newtheorem{lemma}[theorem]{Lemma}
\newtheorem{corollary}[theorem]{Corollary}
\newtheorem{proposition}[theorem]{Proposition}

\theoremstyle{definition}

\newtheorem{assumption}[theorem]{Assumption}
\newtheorem{definition}[theorem]{Definition}

\theoremstyle{remark}
\newtheorem{remark}[theorem]{Remark}

\numberwithin{equation}{section}

\copyrightinfo{2026}{Bobo Hua and Jin Sun}

\newcommand{\partiald}[2]{  \frac{\partial #1 }{\partial #2}}
\newcommand{\partialdd}[2]{ \partial #1/ \partial #2}
\newcommand{\as}[1]{\left\langle #1\right\rangle}
\newcommand{\av}[1]{\left\vert #1\right\vert}

\newcommand{\R}{\mathbb{R}}
\newcommand{\Dom}{\mathrm{Dom}}

\newcommand{\LRc}[1]{\left(#1\right)}

\newcommand{\Om}{\Omega}

\providecommand{\eat}[1]{}

\newcommand{\Hm}[1]{\leavevmode{\marginpar{\tiny%
			$\hbox to 0mm{\hspace*{-0.5mm}$\leftarrow$\hss}%
			\vcenter{\vrule depth 0.1mm height 0.1mm width \the\marginparwidth}%
			\hbox to 0mm{\hss$\rightarrow$\hspace*{-0.5mm}}$\\\relax\raggedright
			#1}}}



\thanks{B. Hua is supported by NSFC, no.~12371056.}

\begin{document}
	
	\title{The hot spots conjecture on Gaussian Spaces}
    
        \author{Bobo Hua}
        \address{Bobo Hua, School of Mathematical Sciences, Fudan University, 200433, Shanghai, China.}
        \email{bobohua@fudan.edu.cn}
        
	\author{Jin Sun}
	\address{Jin Sun, School of Mathematical Sciences, Fudan University, 200433, Shanghai, China.}
	\email{jsun22@m.fudan.edu.cn}

	\begin{abstract}
		We study the hot spots conjecture for domains in the Gaussian space $(\mathbb{R}^n, (2\pi)^{-n/2} e^{-|x|^2/2} dx)$ for $n \ge 2$. Given a bounded domain $\Omega$ with a piecewise smooth boundary, we consider the first nontrivial eigenfunction of the Ornstein--Uhlenbeck operator $L_\gamma = \Delta - \langle x, \nabla \rangle$ subject to Neumann or mixed Dirichlet--Neumann boundary conditions, and prove that its extrema are attained only on the boundary $\partial\Omega$.

        More precisely, we establish the conjecture for two classes of domains: (i) lip domains in Gaussian spaces with mixed boundary conditions, and (ii) $n$-symmetric domains whose intersection with some orthant is a lip domain. As a corollary, we show that any first nontrivial Neumann eigenfunction of a $2$-symmetric domain in the two-dimensional Gaussian space has no interior extrema, provided the second Neumann eigenvalue is simple. Our approach is based on a variational principle for the Hodge Laplacian on weighted manifolds and the Hodge decomposition of differential $1$-forms on Lipschitz domains, extending the variational method of Kennedy--Rohleder \cite{KennedyJamesB2024Oths} from the Euclidean setting to Gaussian spaces. Although de Dios Pont \cite{pont2024convex} has shown that the hot spots conjecture can fail for certain convex domains endowed with suitable log-concave measures, our results identify broad classes of domains for which the conjecture remains valid in Gaussian spaces.
	\end{abstract}

	\maketitle
	
	\section{Introduction}

        The hot spots conjecture, first proposed by J. Rauch in 1974 \cite{rauch1975lecture}, asserts that for any bounded domain $ D \subset \mathbb{R}^n $ with a piecewise smooth boundary, every first nontrivial eigenfunction $\psi_2$ of the Neumann Laplacian attains its maximum and minimum only on the boundary $\partial D$. More precisely, for
        \begin{equation*}
        	\left\{
        	\begin{aligned}
        		-\Delta \psi_2 &= \mu_2\psi_2 \ \ \text{in}\ D,\\
        		\partial_{\mathbf{n}}{\psi_2} &= 0 \quad \text{on}\ \partial D,
        	\end{aligned}\right.
        \end{equation*}
        where $\mu_2$ is the second Neumann eigenvalue, $\mathbf{n}$ denotes the outer unit normal vector and $\Delta:=\sum_{i=1}^{n}\partial^2/\partial x_{i}^2$, then
        \begin{align*}
        	 \max_{x\in\partial D}\psi_2=\max_{x\in\overline{D}}\psi_2. 
        \end{align*}
        
        This conjecture has been the subject of extensive research in spectral geometry. It has been verified for several important classes of domains, including triangles \cites{MR4045963,MR5015432,MR3369351}, lip domains \cite{MR2051611}, convex domains with one axis of symmetry \cite{MR1926894}, and domains with two axes of symmetry \cite{MR1775736}; see also \cites{MR1788041,MR1694534,MR2048182,MR2572703,MR4001630,MR1873295,MR4944076,MR2092873} for further results. On the other hand, counterexamples have been constructed for domains with holes \cites{MR2169871,MR1680567}, and de Dios Pont \cite{pont2024convex} has recently demonstrated failures for certain convex sets in high dimensions. Numerical evidence of failure has been given in \cite{MR4349243}, and sharp quantitative bounds on the extent to which the conjecture can fail were obtained in \cite{de2025sharp}. Recent developments on the analysis of critical points and the location of hot spots have further deepened the understanding of Neumann eigenfunctions; see \cites{MR4462187,MR4107000,MR4615892,MR1937899,MR3022811,MR4625002,MR5030105,DGJYYZ2026}. The conjecture has also been studied under mixed Dirichlet–Neumann boundary conditions, both in the plane and in higher dimensions \cites{MR4860866,MR4855877,MR4949380,KennedyJamesB2024Oths}, as well as in the discrete setting \cites{MR4315494,MR4806996,MR3912970}.

        Despite this extensive body of work in Euclidean spaces, comparatively little is known in the setting of Riemannian manifolds. Krej\v{c}i\v{r}\'i{k}  and Tu\v{s}ek \cite{MR3912674} proved that the hot spots conjecture holds for thin curved strips on surfaces. Freitas \cite{MR1909291} constructed some counterexamples on the unit disk equipped with a certain $S^1$-invariant metric. More recently, Hatcher proved that the hot spots conjecture holds for all non-acute geodesic triangles of constant negative curvature \cite{hatcher2025hot1} and for certain cones and warped product manifolds \cite{hatcher2025hot2}. The general case remains open.
        
        The goal of this paper is to establish the hot spots conjecture in the setting of Gaussian spaces, thereby extending the Euclidean results for lip domains and symmetric domains to a natural class of weighted Riemannian manifolds. Our principal innovation is the development of a variational framework based on differential $1$-forms and the Hodge decomposition on Lipschitz domains in Gaussian spaces. This approach, inspired by the variational method of Kennedy and Rohleder \cites{MR4738461,KennedyJamesB2024Oths} for Euclidean vector fields, yields an intrinsic proof that does not rely on explicit eigenfunction computations or analysis of critical points.
        
        We denote the Gaussian measure on $\mathbb{R}^n$ by
        \begin{equation}\label{E:Gaussian_measure}
            d\gamma = \frac{1}{(2\pi)^{n/2}}e^{-{\av{x}^2}/{2}},
        \end{equation}
        and denote the Gaussian space of dimension $n$ by $\mathcal{M}^n=(\R^n,d\gamma)$. The Ornstein-Uhlenbeck operator is given by
        $$
        L_\gamma u:= e^{|x|^2/2}\operatorname{div}(e^{-|x|^2/2}\nabla u)= \Delta u - \langle x, \nabla u \rangle
        $$
        We consider the following eigenvalue problem of $L_\gamma$ with a mixed boundary condition on a bounded Lipschitz domain $\Omega$, namely,
        \begin{equation}\label{E:OU_equation}
            \begin{cases} 
                -\Delta u(x) + \left\langle x, \nabla u(x) \right\rangle=\lambda u(x), & x \in \Omega, \\ 
                u = 0, & \text{on}\ \Sigma\\
                \partial_\mathbf{n} u = 0 &\text{on}\ \Sigma^\prime,
            \end{cases}
        \end{equation}
         where $\mathbf{n}$ is the outer normal vector, $\partial\Omega = \overline{\Sigma}\cup\overline{\Sigma'}$ and $\Sigma\cap\Sigma'=\emptyset$.

        We say that a bounded domain $\Omega\subset\mathbb{R}^n$ has \emph{piecewise $C^\infty$ boundary} if $\Omega$ is a Lipschitz domain and there exist finitely many pairwise disjoint connected relatively open $C^\infty$ hypersurfaces $P_1,\ldots,P_N\subset\partial\Omega$ such that
        \begin{align*}
            \partial\Omega=\bigcup_{j=1}^N \overline{P}_j,
        \end{align*}
        and the singular set
        \begin{align*}
            \Sigma=\partial\Omega\setminus\bigcup_{j=1}^N P_j
        \end{align*}
        is contained in a finite union of embedded $C^\infty$ submanifolds of dimension at most $n-2$. On each face $P_j$, the outer unit normal and the second fundamental form are defined in the classical sense.
        
        Throughout the paper, we impose the following standing assumption.
        \begin{assumption}\label{assumption}
        	The domain $\Omega\subset\R^n$ is a bounded Lipschitz domain that satisfies a uniform exterior ball condition and has piecewise $C^\infty$ boundary.
        \end{assumption}
        Under the uniform exterior ball condition, the domains of the Dirichlet Laplacian and the Neumann Laplacian on $L^2(\Omega)$ are embedded in the Sobolev space $H^2(\Omega)$; see \cite{MR2581375}.

        Let $e_j:=\partialdd{}{x_j}\in T\mathcal{M}^n$. Then for a $1$-form $u\in \Lambda^1T\mathcal{M}^n$, each component of $u$ is given by $u_j:=u(e_j)$ for each $j=1,\ldots,n$. Hence, every $1$-form $u$ can be written by
        \begin{align}\label{expression}
        	u = \sum_{j=1}^{n}u_jd x_j.
        \end{align} 
        
        We now introduce the key geometric notions underlying our results.
        
        \begin{definition}[Lip domain]\label{D1}
            A Lipschitz domain $ \Omega\subset \mathcal{M}^n$ is called a \textit{lip domain} if after a suitable orthogonal transformation (rotation or reflection), its boundary admits a decomposition $\partial\Omega=\overline{\Gamma}_1\cup\overline{\Gamma}_2$, where $\Gamma_1$ and $\Gamma_2$ are  open subsets defined as follows. Let $\mathbf{n}(x)$ denote the outward unit normal vector to $\Omega$ at $x$, and set
            \begin{align*}
            	S_1&:=\{x\in \partial\Omega:\mathbf{n}(x)\ \text{has exactly two nonzero components with opposite signs}\},\\
            	S_2&:=\{x\in \partial\Omega:x_j=0\ \text{and}\ \mathbf{n}(x)= \pm e_j \ \text{for some}\ j\}.
            \end{align*}        
            Then
            \begin{align*}
            	\Gamma_1:=\mathrm{int}\LRc{\overline{S_1}},\ 
            	\Gamma_2:=\mathrm{int}\LRc{\overline{S_2}}.
            \end{align*}
        \end{definition}   

        \begin{figure}[htbp]
            \centering
            \begin{subfigure}[b]{0.45\textwidth}
                \centering
                \includegraphics[width=\textwidth]{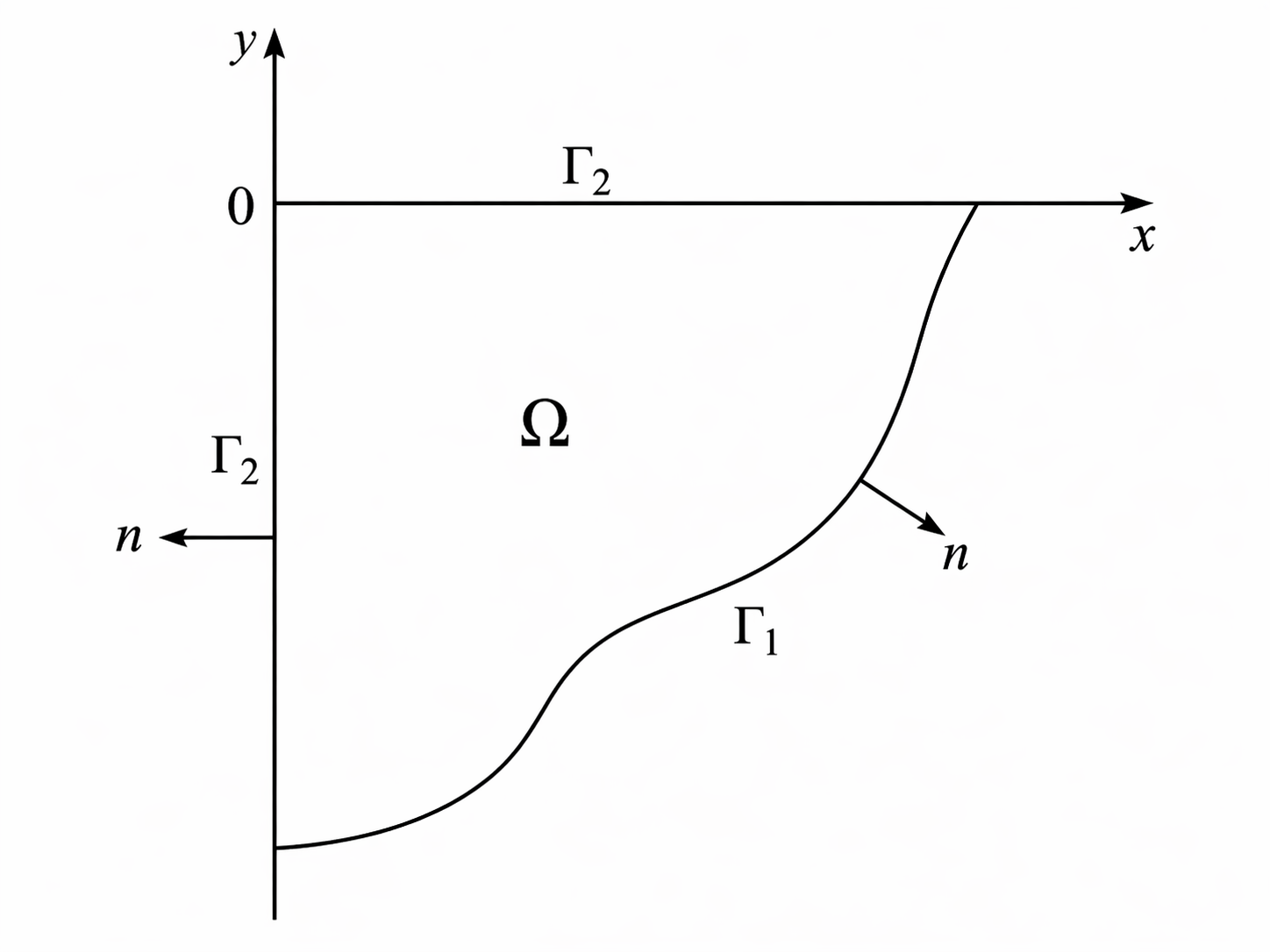}
                \caption{Lip domain in $\mathcal{M}^2$}
                \label{fig:lip2d}
            \end{subfigure}
            \hfill
            \begin{subfigure}[b]{0.45\textwidth}
                \centering
                \includegraphics[width=\textwidth]{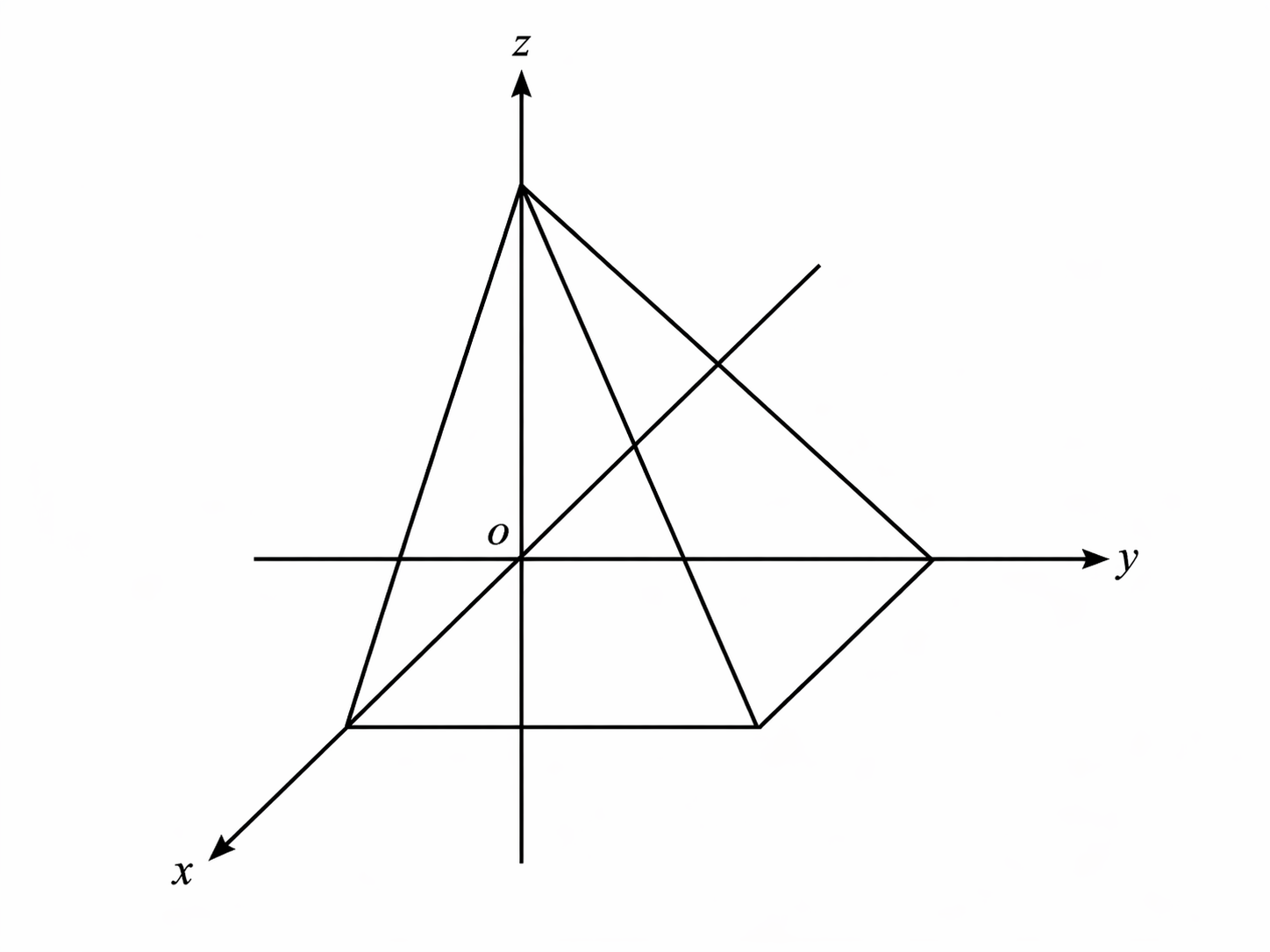}
                \caption{Lip domain in $\mathcal{M}^3$}
                \label{fig:lip3d}
            \end{subfigure}
            \caption{Examples of lip domains}
            \label{fig:lipdomains}
        \end{figure}
        
        \begin{definition}[$n$-symmetric domain]\label{D2}
        	A Lipschitz domain $ \Omega \subset \mathcal{M}^n$ is \textit{$n$-symmetric} if $\Omega$ is invariant under reflection in each coordinate hyperplane $ \{x_k = 0\}, k = 1,\ldots, n $.
        	Given an $n$-symmetric domain $\Omega$, we say that a subset $\mathcal{O}$ is an \textit{orthant of $\Omega$} if $\mathcal{O}=\{x\in \Omega:(-1)^{i_1}x_1>0,\ldots,(-1)^{i_n}x_n>0\}$ for some integers $i_1,\ldots,i_n$.
        \end{definition}    

        \begin{figure}[htbp]
            \centering
            \begin{subfigure}[b]{0.45\textwidth}
                \centering
                \includegraphics[width=\textwidth]{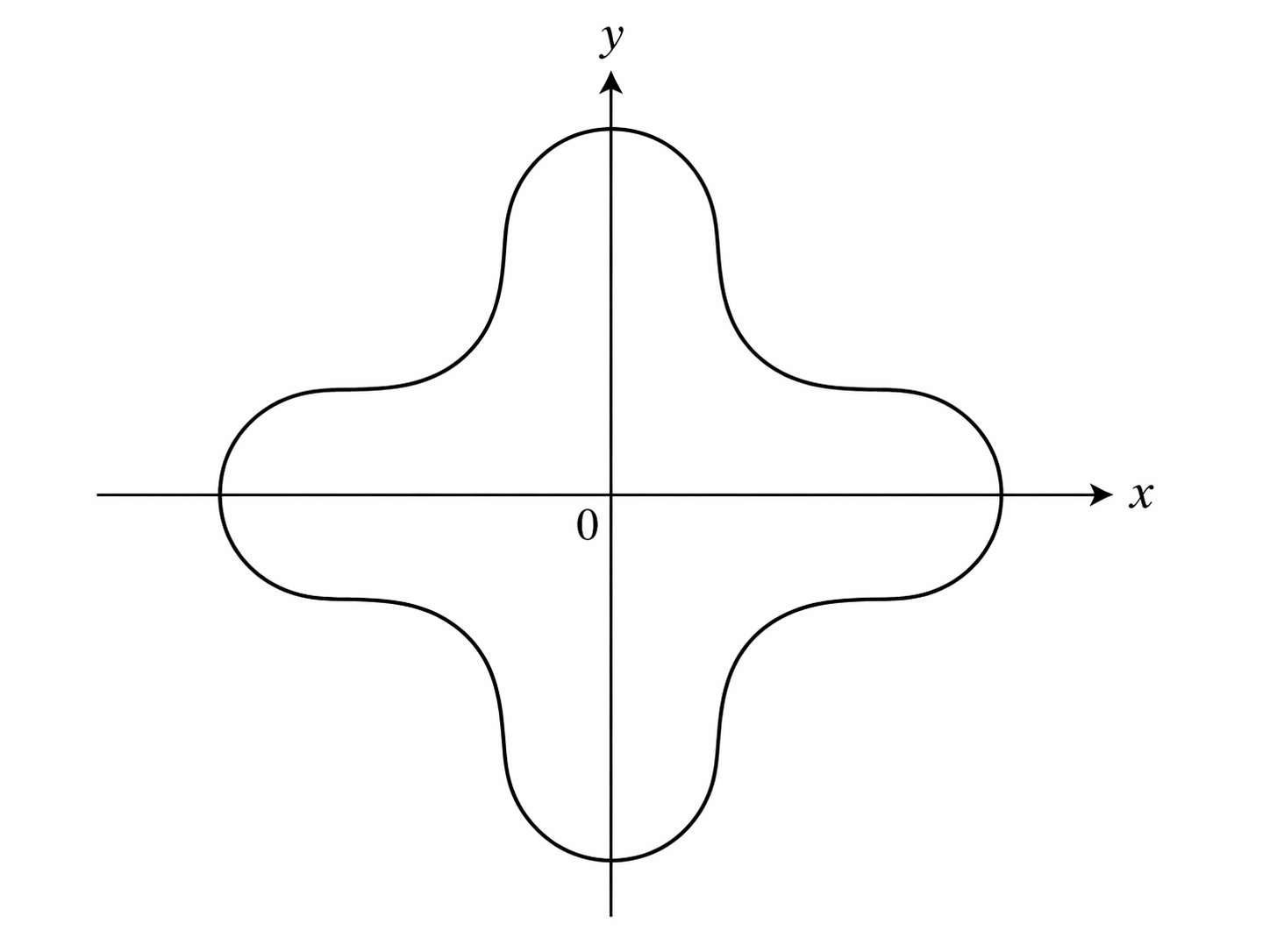}
                \caption{$2$-symmetric domain in $\mathcal{M}^2$}
                \label{fig:symmetric2d}
            \end{subfigure}
            \hfill
            \begin{subfigure}[b]{0.35\textwidth}
                \centering
                \includegraphics[width=\textwidth]{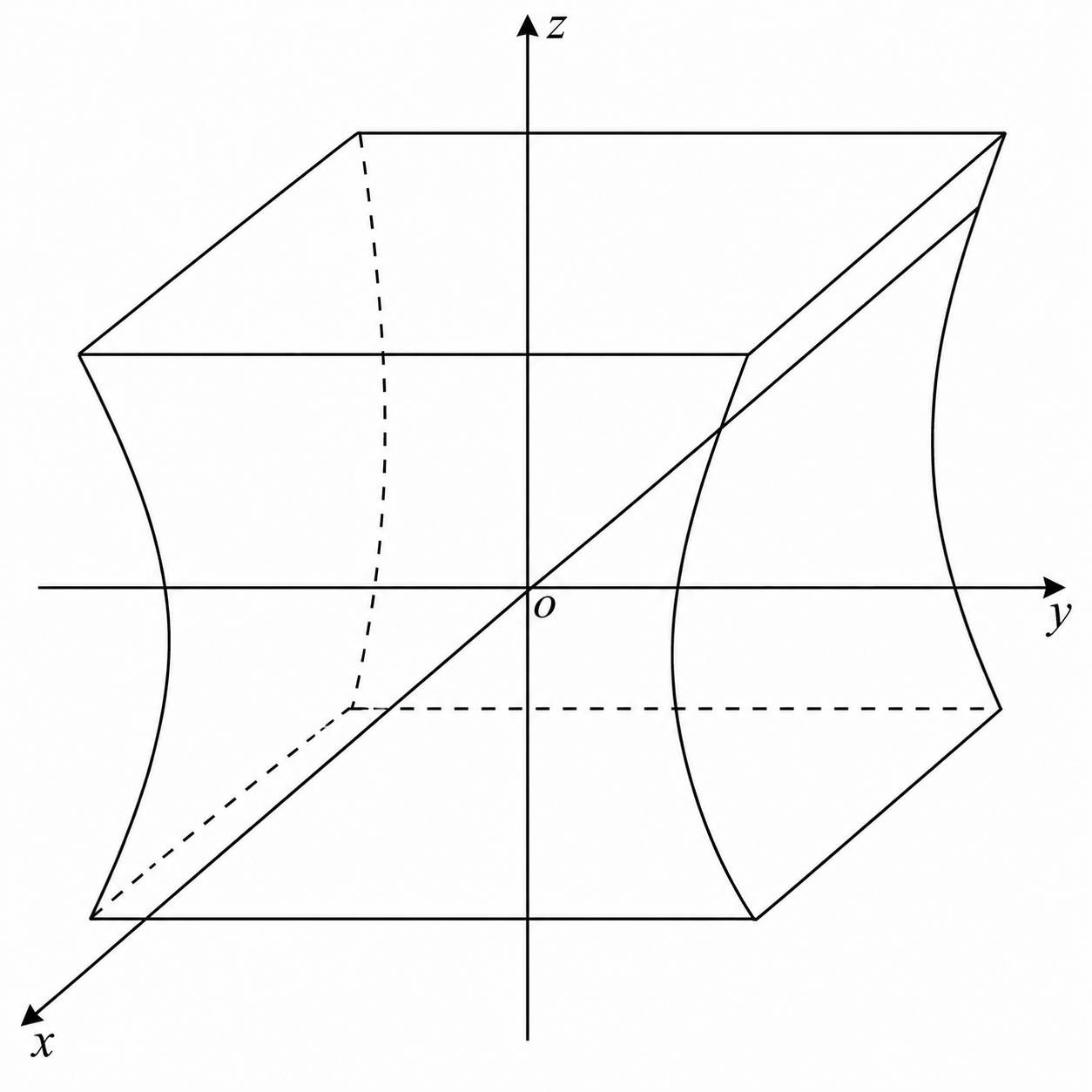}
                \caption{$3$-symmetric domain in $\mathcal{M}^3$}
                \label{fig:symmetric3d}
            \end{subfigure}
            \caption{Examples of $n$-symmetric domains}
            \label{fig:symmetricdomains}
        \end{figure}
        
        \begin{remark}
            Due to the drift term $-x\cdot\nabla$ in the Ornstein–Uhlenbeck operator, the symmetry axes are required to pass through the origin. This differs slightly from the Euclidean setting treated in \cite{KennedyJamesB2024Oths}.
        \end{remark}

        \begin{definition}[Antisymmetric eigenfunction]
            Let $ \Omega \subset \mathcal{M}^n$ be an $n$-symmetric Lipschitz domain. An eigenfunction $w$ of problem \eqref{E:OU_equation} with $\Sigma=\emptyset$ (pure Neumann conditions) is called \textit{$k$-antisymmetric} if $w$ is odd with respect to the hyperplane $\{x_k=0\}$. A $k$-antisymmetric eigenfunction $w_k$ is called a \textit{first $k$-antisymmetric eigenfunction}  if it corresponds to the smallest eigenvalue among all $k$-antisymmetric eigenfunctions.
        \end{definition}

        As in \cite{MR2839867}, we introduce the class of creased domains.
        \begin{definition}[Creased domain]\label{D:creased_domain}
            Let $\Omega \subset \mathcal{M}^n$ be a bounded Lipschitz domain. Suppose that two disjoint open subsets $\Sigma$ and $\Sigma'$ of the boundary $\partial\Omega$ (endowed with the relative topology induced by $\mathbb{R}^n$) are given such that  $\partial\Omega = \overline{\Sigma} \cup \overline{\Sigma'}$. Take any point $x_o \in \overline{\Sigma}\cap\overline{\Sigma'}$. We say that $\Omega$ is \emph{creased near $x_o$} if there exist constants $\kappa>0$ and $r>0$ for which the following conditions are satisfied. There exist two Lipschitz functions  
            \begin{align*}
                \varphi:\mathbb{R}^{n-1}\to\mathbb{R}, \qquad \xi:\mathbb{R}^{n-2}\to\mathbb{R}
            \end{align*}
            such that, after a suitable rotation and translation of coordinates (implicitly assumed in the representation below), we have  
            \begin{align*}
                \Omega\cap B_{2r}(x_o) &= B_r(x_o)\cap\Bigl\{ (x_1,x'',x_n)\in\mathbb{R}\times\mathbb{R}^{n-2}\times\mathbb{R} \;:\; x_n > \varphi(x_1,x'') \Bigr\},\\[4pt]
                \Sigma\cap B_{2r}(x_o) &= B_r(x_o)\cap\partial\Omega\cap\Bigl\{ (x_1,x'',x_n) \;:\; x_1 < \xi(x'') \Bigr\},\\[4pt]
                \Sigma'\cap B_{2r}(x_o) &= B_r(x_o)\cap\partial\Omega\cap\Bigl\{ (x_1,x'',x_n) \;:\; x_1 > \xi(x'') \Bigr\},
            \end{align*}
            and moreover the partial derivative $\partial_{x_1}\varphi$ satisfies  
            \begin{align*}
                \partial_{x_1}\varphi &\ge \kappa \quad \text{on } \bigl\{ (x_1,x'')\in\mathbb{R}\times\mathbb{R}^{n-2} \;:\; x_1 < \xi(x'') \bigr\},\\[4pt]
                \partial_{x_1}\varphi &\le -\kappa \quad \text{on } \bigl\{ (x_1,x'')\in\mathbb{R}\times\mathbb{R}^{n-2} \;:\; x_1 > \xi(x'') \bigr\}.
            \end{align*}
            Finally, the pair $(\Sigma,\Sigma')$ is called \emph{admissible decomposition of $\partial\Omega$} and $\Omega$ is called \emph{creased relative to $(\Sigma,\Sigma')$} if $\Omega$ is creased near every point $x_o\in\overline{\Sigma}\cap\overline{\Sigma'}$. We simply say that  $\Omega$ is \emph{creased relative to $\Sigma$} if  $\Omega$ is creased relative to $(\Sigma,\partial\Omega\backslash\overline{\Sigma})$.
        \end{definition}
        Informally speaking, the boundary is partitioned into admissible patches $\Sigma$ and $\Sigma'$ that meet at angles $<\pi$.

        We can now state our main results. The first theorem establishes the hot spots conjecture for lip domains in Gaussian spaces under mixed boundary conditions.
        \begin{theorem}\label{Tlip}
            Let $ \Omega\subset \mathcal{M}^n$ be a lip domain satisfying Assumption \ref{assumption}, and let $\Sigma=\mathrm{int}\LRc{{\Gamma}_2\cap\{x_k=0\}}$ be a connected and open subset for some integer $1\leq k\leq n$. Assume that $\Omega$ is creased relative to $\Sigma$. Let $ \psi $ be any eigenfunction of problem \eqref{E:OU_equation} corresponding to the first eigenvalue $\lambda_1^\Sigma$. Then $ \psi $ is monotonic in each coordinate direction $e_i$, $i=1,\ldots,n$. In particular, $ \psi $ attains its maximum and minimum only on $ \partial \Omega $.
        \end{theorem}

        The second theorem treats $n$-symmetric domains and generalizes the Euclidean result of Jerison–Nadirashvili \cite{MR1775736}*{Theorem 1.1}.

        \begin{theorem}\label{Tsymmetric}
            Let $\Omega\subset \mathcal{M}^n$ be an $n$-symmetric domain satisfying Assumption \ref{assumption}. Assume that for some orthant $\mathcal{O}$ of $\Omega$, the intersection $ \Omega\cap\mathcal{O}$ is a lip domain and creased relative to $\partial( \Omega\cap\mathcal{O})\cap\{x_k=0\}$ for each $ k = 1,\ldots,n$. Then for each $k$, the first $k$-antisymmetric eigenfunction $w_k$ is unique up to scalar multiples and is even with respect to each coordinate hyperplane $ \{x_{k^\prime} = 0\}$ for $k^\prime \neq k $. Moreover, $\partial_kw_k\geq 0$ in $\Omega$, and the maximum and minimum of $ w_k $ are attained only on $\partial \Omega $.
        \end{theorem}
        
        When $n=2$, the multiplicity of the second Neumann eigenvalue $\mu_2$ is at most $2$. If $\mu_2$ is simple, then by the symmetry of the domain, any first nontrivial Neumann eigenfunction is necessarily $k$-antisymmetric for some $k\in\{1,2\}$. As a direct consequence of Theorem \ref{Tsymmetric}, we have the following corollary.
        \begin{corollary}\label{Corollary}
            Let $\Omega\subset \mathcal{M}^2$ be a bounded convex domain with piecewise smooth boundary, satisfying the uniform exterior ball condition, and symmetric with respect to both coordinate axes. If the second Neumann eigenvalue $\mu_2(\Omega)$ is simple, then every first nontrivial Neumann eigenfunction attains its extrema only on $\partial\Omega$.
        \end{corollary}
        
        This paper is organized as follows. In Section~\ref{sec2}, we give the Hodge decomposition on Lipschitz domains in Gaussian spaces and establish the function spaces needed for our analysis. Section~\ref{sec3} develops the variational principles for the Gaussian mixed Hodge Laplacian, including the min-max characterization of its eigenvalues. In Section~\ref{sec4}, we apply these principles to lip domains: we prove the key inequality $\eta(\av{u})\leq \eta(u)$ for the Rayleigh quotient of $1$-forms and use it to establish the monotonicity of eigenfunctions, thereby proving Theorem \ref{Tlip}. Finally, Section~\ref{sec5} reduces the case of $n$-symmetric domains to lip domains via a reflection argument, completing the proofs of Theorem \ref{Tsymmetric} and Corollary \ref{Corollary}. We conclude with Proposition \ref{prop}, which shows that the variational method does not extend to lip domains in spaces of constant curvature directly.

        \section{Hodge Decomposition in Gaussian Spaces}\label{sec2}
        Let $\Omega\subset\R^n$ be a lip domain satisfying Assumption \ref{assumption}. Let $(\Sigma,\Sigma')$ be admissible decomposition given in \ref{D:creased_domain} and let $\Omega$ be creased relative to $(\Sigma,\Sigma')$.
        
        We denote by $\Lambda^k=\Lambda^k(T^*\mathcal{M}^n)$ the bundle of exterior $k$-forms and
        $$\Lambda=\bigoplus_{k=0}^{n}\Lambda^k.$$
        
        We introduce the Hodge star operator as a linear map

        $$
        * : \Lambda^k \longrightarrow \Lambda^{n-k},
        $$
        defined pointwise by the following condition: for any $k$-form $\alpha$ and any $(n-k)$-form $\beta$ on a tangent space $T_p\R^n$,$$
        \alpha \wedge * \beta = \langle \alpha, \beta \rangle_{Euc} \, d{vol}_{Euc} ,
        $$
        where $\langle \cdot, \cdot \rangle_{Euc}$ is the inner product on $\Lambda^k$.
        
        For $u\in\Lambda^l$ and $v\in \Lambda^r$ with $r\leq l$, we further define the interior product by
        \begin{align*}
        	v\vee u:=(-1)^{(l-r)(n-l)}*(v\wedge*u).
        \end{align*}
        Then for $w\in\Lambda^{l-r}$, we have the adjoint identity
        \begin{align*}
            \langle v\wedge w,u\rangle=\langle w,v\vee u\rangle .
        \end{align*}
        
        Let $d\colon C^\infty(\Omega,\Lambda^k)\to C^\infty(\Omega,\Lambda^{k+1})$ be the exterior derivative and let $\delta$ be the standard codifferential $\delta\colon C^\infty(\Omega,\Lambda^{k+1})\to C^\infty(\Omega,\Lambda^k)$, given by
        \begin{align*}
        	\delta w = (-1)^{nk+1}*d(*w)\,,
            \qquad w\in C^\infty(\Omega,\Lambda^{k+1}).
        \end{align*}
        
        We define the \emph{Gaussian codifferential} by
        \begin{equation}\label{eq:delta_phi}
            d^* w := e^{{\av{x}^2}/{2}}\delta (e^{-{\av{x}^2}/{2}}w)\,,
            \qquad w\in C^\infty(\Omega,\Lambda^{k+1}).
        \end{equation}
        Therefore, $(d^*)^2 = 0$. For a $k+1$-form $w\in\Lambda^{k+1}$ and a $k$-form $v\in\Lambda^{k}$ with compact supports in $\Om$,
        \begin{align*}
        	(w,dv)_\Omega = \int_{\Omega}\as{w,dv}_{\text{Euc}}d\gamma = \int_{\Omega}\as{d^*w,v}_{\text{Euc}}d\gamma = (d^*w,v)_\Omega.
        \end{align*}
        
        For admissible patches $\Sigma$ and $\Sigma^\prime$ such that $
        \overline{\Sigma}\cup\overline{\Sigma^\prime}=\partial\Omega$, we define
        \begin{align*}
            &L^2_\gamma(\Omega,\Lambda):=  \left\{ u : \Omega \to \Lambda \;\Bigg|\; \int_{\Omega} \langle u, u \rangle_{\text{Euc}} \, d\gamma < \infty \right\},\\
        	&H_\gamma(\Omega,\Sigma,d):=\{u\in L^2_\gamma(\Omega,\Lambda):d u\in L^2_\gamma(\Omega,\Lambda), u_{tan}=0\ \text{on}\ \Sigma\},\\
        	&H_\gamma(\Omega,\Sigma^\prime,d^*):=\{u\in L^2_\gamma(\Omega,\Lambda):d^* u\in L^2_\gamma(\Omega,\Lambda), u_{nor}=0\ \text{on}\ \Sigma^\prime\},\\
        	&\mathscr{H}_\gamma(\Omega, \Sigma, \Sigma'):=\{u\in L^2_\gamma(\Omega,\Lambda):d u=0, d^* u=0\ \text{in}\ \Omega, u_{tan}\big|_{\Sigma}=0, u_{nor}\big|_{\Sigma^\prime}=0\,\},
        \end{align*}
        where $\nu$ is the outward unit conormal $1$-form on $\partial\Omega$ and
        \begin{align*}
        	u_{tan}:=\nu\vee(\nu\wedge u),\  u_{nor}:=\nu\wedge(\nu\vee u).
        \end{align*}
        By definition, $u=u_{tan}+u_{nor}$. 
        
        The spaces $H_\gamma(\Omega,\Sigma,d)$ and $H_\gamma(\Omega,\Sigma^\prime,d^*)$ are Hilbert spaces under the graph norms
        \begin{align*}
            \big\|u\big\|_{H_\gamma(\Omega,\Sigma,d)}^2 &:= \big\|u\big\|_{L^2_\gamma(\Omega)}^2+\big\|du\big\|_{L^2_\gamma(\Omega)}^2,
            \\
            \big\|u\big\|_{H_\gamma(\Omega,\Sigma^\prime,d^*)}^2 &:= \big\|u\big\|_{L^2_\gamma(\Omega)}^2+\big\|d^*u\big\|_{L^2_\gamma(\Omega)}^2.
        \end{align*}
  
        The boundary conditions in these spaces are understood in the sense of traces of the appropriate Sobolev spaces on Lipschitz domains, as developed in \cite{MR2839867}. 
  
        On a bounded domain $\Omega$, the weighted space $L^2_\gamma(\Omega,\Lambda^k)$ is norm-equivalent to the standard space $L^2(\Omega,dx,\Lambda^k)$. Consequently, the natural inclusion
        \begin{align*}
            H_\gamma(\Omega, \Sigma, d) \cap H_\gamma(\Omega,\Sigma^\prime,d^*)\;\hookrightarrow\; L^2_\gamma(\Omega,\Lambda)
        \end{align*}
        is compact, where the former space is equipped with the norm
        \begin{align}\label{E:norm_H_1}
            u\mapsto \big\|u\big\|_{L^2_\gamma(\Omega)}+\big\|du\big\|_{L^2_\gamma(\Omega)}+\big\|d^*u\big\|_{L^2_\gamma(\Omega)}.
        \end{align}
        This compactness result can be derived from the Gaffney’s inequality (see, e.g., \cite{MR1367287}).
        
        Following \cite{MR2839867}*{Proposition 2.9}, we obtain the following Hodge decomposition in Gaussian spaces.
       \begin{theorem}\label{T:Hodge_decomposition}
           Let $\Omega\subset\R^n$ be a Lipschitz domain and $(\Sigma,\Sigma')$ be admissible decomposition of $\partial\Omega = \overline{\Sigma}\cup\overline{\Sigma'}$. Then
           \begin{enumerate}[label=\textup{(\roman*)},leftmargin=2em]
                \item The spaces $dH_\gamma(\Omega, \Sigma, d)$ and $d^* H_\gamma(\Omega, \Sigma', d^*)$ are closed in $L^2_\gamma(\Omega, \Lambda)$, and the space $\mathscr{H}^k_\gamma(\Omega,\Sigma,\Sigma')$ is finite-dimensional.
                \item The following Hodge decomposition holds:
                \begin{equation}\label{Hodgedecomposition}
        	           L^2_\gamma(\Omega, \Lambda) = dH_\gamma(\Omega, \Sigma, d) \oplus d^* H_\gamma(\Omega, \Sigma', d^*) \oplus \mathscr{H}_\gamma(\Omega, \Sigma, \Sigma'),
                \end{equation}
            where the direct sums are orthogonal.
            \end{enumerate}
          \end{theorem}
   
        \section{Variational Principles in Gaussian Spaces}\label{sec3}
        
        We define unbounded linear operators $d_\Sigma:L^2_\gamma(\Omega, \Lambda)\rightarrow L^2_\gamma(\Omega, \Lambda)$ and $d^*_{\Sigma^\prime}:L^2_\gamma(\Omega, \Lambda)\rightarrow L^2_\gamma(\Omega, \Lambda)$ by 
        \begin{align*}
        	d_\Sigma&:=d u,\quad \forall u\in \mathrm{Dom}(d_\Sigma)=H_\gamma(\Omega,\Sigma,d);\\
        	d^*_{\Sigma^\prime}&:=d^* u,\quad \forall u\in \mathrm{Dom}(d^*_{\Sigma^\prime})=H_\gamma(\Omega,\Sigma^\prime,d^*).
        \end{align*}
        The direct calculation yields $(d_\Sigma)^*=d^*_{\Sigma^\prime}$. 
        
        Since the norm on $L^2_\gamma(\Omega,\Lambda^k)$ is equivalent to the norm on $L^2(\Omega,dx,\Lambda^k)$ for any bounded domain $\Omega$, from the abstract theory \cite{MR2839867}*{Section 2}, we have the following theorem.
        \begin{theorem}[Gaussian analogue of {\cite{MR2839867}*{Theorem~4.5}}]\label{thm:Hodge_Laplace}
            Let $\Omega\subset\R^n$ be a bounded Lipschitz domain with admissible decomposition $\partial\Omega=\overline\Sigma\cup\overline{\Sigma'}$.  Then the following hold.

            \begin{enumerate}[label=\textup{(\roman*)},leftmargin=2em]
                \item\label{it:sa} The \emph{Gaussian mixed Hodge Laplacian}, defined in the sense of composition of unbounded operators by
                \begin{equation}\label{eq:HodgeLap}
                    \Delta_{\gamma,\Sigma}:=d_{\Sigma} \circ d^*_{\Sigma^\prime} + d^*_{\Sigma^\prime} \circ d_{\Sigma},
                \end{equation}
                is self-adjoint and nonnegative on $L^2_\gamma(\Omega, \Lambda)$ and has a compact resolvent.

                \item\label{it:dense} Its domain $\Dom(\Delta_{\gamma,\Sigma})$ is dense in $H_\gamma(\Omega, \Sigma, d) \cap H_\gamma(\Omega,\Sigma^\prime,d^*)$ when the latter is equipped with the norm \eqref{E:norm_H_1}.

                \item\label{it:form} $\Dom\bigl((\Delta_{\gamma,\Sigma})^{1/2}\bigr) = H_\gamma(\Omega, \Sigma, d) \cap H_\gamma(\Omega,\Sigma^\prime,d^*)$.

                \item\label{it:sectorial} $\Delta_{\gamma,\Sigma}$ is sectorial and $-\Delta_{\gamma,\Sigma}$ generates an analytic semigroup on $L^2_\gamma(\Omega, \Lambda)$.

                \item\label{it:spectrum} Consider the eigenform of Gaussian mixed Hodge Laplacian problem with mixed boundary condition as follows:
                \begin{equation}\label{$1$-form}
        	        \Delta_{\gamma,\Sigma} u = \eta u,\qquad u\in H_\gamma(\Omega, \Sigma, d) \cap H_\gamma(\Omega,\Sigma^\prime,d^*).
                \end{equation}
                The spectrum of $\Delta_{\gamma,\Sigma}$ is discrete and consists of only nonnegative real eigenvalues with no finite accumulation point.  When listed non-decreasingly according to multiplicity, these satisfy the min-max principle
                \begin{equation*}
        	        \eta_j = \min_{\substack{L_j \text{ subspace of } H_\gamma(\Omega, \Sigma, d) \cap H_\gamma(\Omega,\Sigma^\prime,d^*) \\ \dim(L_j) = j}} \left( \max_{0 \neq u \in L_j} \eta(u) \right), \quad j \in \mathbb{N},
                \end{equation*}
                where
                \begin{equation}\label{variation}
        	       \eta(u):=\frac{\|du\|^2_{L^2_\gamma(\Omega, \Lambda)} + \|d^* u\|^2_{L^2_\gamma(\Omega, \Lambda)}}{\|u\|^2_{L^2_\gamma(\Omega, \Lambda)}}.
                \end{equation}
            \end{enumerate}
        \end{theorem}

        \begin{remark}
            A direct calculation shows that for any function $f\in C^{\infty}(\Omega)\cap H_\gamma(\Omega, \Sigma, d) \cap H_\gamma(\Omega,\Sigma^\prime,d^*)$,
            \begin{align*}
                \Delta_{\gamma,\Sigma} f= -L_\gamma f.
            \end{align*}
            Therefore, to study the eigenfunctions of problem \eqref{E:OU_equation}, we only need to consider the eigenforms of problem \eqref{$1$-form}.
        \end{remark}
        
        In the case $n=2$, the space of Gaussian harmonic $1$-forms on any simply connected domain is trivial, so that we have $\mathscr{H}_\gamma(\Omega, \Sigma, \Sigma')=\{0\}$ for any bounded and simply connected domain $\Omega\subset\mathcal{M}^n$ with $\Sigma\subset\partial\Omega$ connected. Let $\psi_k^\Sigma$ be the $k$-th eigenfunction with Dirichlet condition on $\Sigma$ and Neumann condition on $\Sigma^\prime$ associated with eigenvalue $\lambda_k^\Sigma$, and let $\varphi_k^\Sigma$ be the $k$-th eigenfunction with Dirichlet condition on $\Sigma^\prime$ and Neumann condition on $\Sigma$ associated with eigenvalue $\lambda_k^{\Sigma^\prime}$. That is,
        \begin{equation}\label{mixedequation}
        	\left\{
        	\begin{aligned}
        		\Delta_{\gamma,\Sigma}  \psi_k^\Sigma&= \lambda_k^{\Sigma} \psi_k^\Sigma\quad \text{in}\ \Omega,\\
        		{\psi_k^\Sigma} &= 0 \quad\text{on}\ \Sigma,\\
        		\partial_{\mathbf{n}}{\psi_k^\Sigma} &= 0 \quad\text{on}\ \Sigma^\prime,
        	\end{aligned}\right.
        \end{equation}
        and 
        \begin{equation}\label{mixedequation2}
        	\left\{
        	\begin{aligned}
        		\Delta_{\gamma,\Sigma}  \varphi_k^\Sigma&= \lambda_k^{\Sigma^\prime} \varphi_k^\Sigma\quad \text{in}\ \Omega,\\
        		{\varphi_k^\Sigma} &= 0 \quad\text{on}\ \Sigma^\prime,\\
        		\partial_{\mathbf{n}}{\varphi_k^\Sigma} &= 0 \quad\text{on}\ \Sigma.
        	\end{aligned}\right.
        \end{equation}
        
        By the spectral theorem and by the commutativity $\Delta_\gamma \circ d\, f = d \circ \Delta_\gamma f$ and $\Delta_\gamma \circ d^* (fd vol_{\mathcal{M}^2})= d^* \circ \Delta_\gamma(fd vol_{\mathcal{M}^2})$ for $f\in C^{\infty}(\Omega)$, we deduce that $1$-forms $\{d\psi_k^\Sigma\}_{k=1}^\infty$ and $1$-forms $\{d^*(\varphi_k^\Sigma d vol_{\mathcal{M}^2})\}_{k=1}^\infty$ form orthonormal eigenbases of $dH_\gamma(\Omega, \Sigma, d)$ and $d^* H_\gamma(\Omega,\Sigma^\prime,d^*)$ respectively, hence they form an orthonormal eigenbasis of $L^2_\gamma(\Omega, \Lambda^1)$. 
        Therefore, we obtain
        \begin{align*}
        	\left\{\eta_j\right\}_{j=1}^\infty=\{\lambda_k^{\Sigma}\}_{k=1}^\infty\cup\{\lambda_k^{\Sigma^\prime}\}_{k=1}^\infty.
        \end{align*}
        
        In the case $n>2$, $d^* H(\Omega,\Sigma^\prime,d^*)$ is complicated and cannot be expressed explicitly. However, the commutativity $\Delta_\gamma \circ d\, f = d \circ \Delta_\gamma f$ still holds for any $0$-form $f\in C^{\infty}(\Omega)$. Hence, the $1$-forms $\{d\psi_k^\Sigma\}_{k=1}^\infty$ always form an orthonormal eigenbasis of $dH_\gamma(\Omega, \Sigma, d)$, and
        \begin{align*}
        	\{\lambda_k^{\Sigma}\}_{k=1}^\infty\subset\{\eta_j\}_{j=1}^\infty.
        \end{align*} 
        
        Note that $\lambda_k^{\Sigma}=\mu_k$ and $\lambda_k^{\Sigma^\prime}=\lambda_k$ if $\Sigma=\emptyset$, where $\mu_k$ is the $k$-th Neumann eigenvalue and $\lambda_k$ is the $k$-th Dirichlet eigenvalue of $\Omega$.

        \section{Lip Domains Have No Hot Spots}\label{sec4}

        In this section, we discuss the hot spots of lip domains in $\mathcal{M}^n$. Recall that $u_j:=u(e_j)$ with the standard basis $e_1,\dots,e_n$ of $\mathbb{R}^n$ and $u = \sum_{j=1}^{n}u_jd x_j$ for any $u\in \Lambda^1$. 
        Let the $1$-form $\av{u}$ be defined as
        \begin{align*}
        	\av{u}:=\sum_{j=1}^{n}\av{u_j} d x_j.
        \end{align*}
        We write $u\geq 0$ to mean $u_j\geq 0$ for all $j=1,\ldots,n$. Obviously, $\av{u}\geq 0$ for every $1$-form $u\in \Lambda^1$.

        We say that a point $p\in\partial\Omega$ is regular if $\partial\Omega$ is smooth at $p$. In particular, we say that $\Sigma$ is regular if every point in $\Sigma$ is regular. The following lemma can be found in \cite{KennedyJamesB2024Oths}*{Lemma 4.2}, which asserts that the outer unit conormal has at most two nonzero components in a neighborhood of each regular point $p_0$.
        \begin{lemma}\label{regularpoints}
        	Let $ \Omega\subset \mathcal{M}^n$ be a lip domain satisfying Assumption \ref{assumption}. Then for each regular point $p_0\in\partial\Omega$, there is an open neighborhood $\mathcal{O}(p_0)$ in $\partial\Omega$ and $k,l\in\{1,\ldots,n\}$ such that $\nu_j(p)=0$ for all $p\in\mathcal{O}(p_0)$ and each $j\neq k,l$.
        \end{lemma}
        
        A direct calculation shows that for $n\geq 2$ and $1$-form $u = \sum_{j=1}^{n}u_jd x_j$,
        \begin{align*}
            (\Delta_\gamma u)_i = \Delta_\gamma u_i + u_i = -L_\gamma u_i + u_i.
        \end{align*}
        
        Therefore, inspired by the proof of Theorem 1.7 in \cite{KennedyJamesB2024Oths}, we have the following essential lemma.
        \begin{lemma}\label{L:essential_lemma}
            Let $ \Omega\subset \mathcal{M}^n$ be a lip domain satisfying Assumption \ref{assumption}, and let $\Sigma=\mathrm{int}\LRc{{\Gamma}_2\cap\{x_k=0\}}$ be a connected and open subset for some integer $1\leq k\leq n$. Assume that $\Omega$ is creased relative to $\Sigma$. Let $\mathbf{n}$ be the outer unit normal vector and $\nu$ be the outer unit conormal. Let $\mathrm{II}$ be the second fundamental form on $\partial\Om$, given by $\mathrm{II}(X, Y) := g_{Euc}(\nabla_X \mathbf{n},Y)$ for $X, Y\in T\partial\Omega$. Then for any $1$-form $u\in H_\gamma(\Omega, \Sigma, d) \cap H_\gamma(\Omega,\Sigma^\prime,d^*)$, 
            \begin{equation}\label{E:eta_u_equivalent}
                \eta(u) = \frac{\int_{\Omega} \big( |\nabla \mathbf{u}|^2 + |\mathbf{u}|^2 \big)\,d\gamma + \int_{\partial\Omega}\mathrm{II}(\mathbf{u}^\top, \mathbf{u}^\top)\,d\gamma_\sigma}{\int_{\Omega} |\mathbf{u}|^2 \,d\gamma},
            \end{equation}
            where $\mathbf{u}=\sum_i u_ie_i$, $\mathbf{u}^\top = \mathbf{u} - (\mathbf{u}\cdot\mathbf{n})\mathbf{n}$ is the tangential component of  $\mathbf{u}$, and $d\gamma_{\sigma}=\frac{1}{(2\pi)^{n/2}}e^{-\av{x}^2/2}d\mathcal{H}^{n-1}$ is the weighted surface measure.
            
            In particular, $\av{u}\in H_\gamma(\Omega, \Sigma, d) \cap H_\gamma(\Omega,\Sigma^\prime,d^*)$ and
            \begin{equation}\label{E:absolute_minimizer}
                \eta(\av{u})\leq \eta(u)
            \end{equation}
        \end{lemma}
        \begin{proof}
            For the term $du$, we have
            \begin{align*}
                |du|^2 &= \frac{1}{2}\sum_{i,j=1}^n (\partial_i u_j - \partial_j u_i)^2 \\
                &= \sum_{i,j} (\partial_i u_j)^2 - \sum_{i,j} \partial_i u_j \partial_j u_i.
            \end{align*}
            The first sum is exactly $|\nabla \mathbf{u}|^2 = \sum_{i,j} (\partial_i u_j)^2$.

            For $ d^*u$ we have
            $$
            (d^* u)^2 = (\operatorname{div}\mathbf{u})^2 - 2(\operatorname{div}\mathbf{u})(x\cdot \mathbf{u}) + (x\cdot \mathbf{u})^2,
            $$
            where $\operatorname{div}\mathbf{u} = \sum_i \partial_i u_i$.

            Adding the two expressions gives
            \begin{align}
                \|du\|_{L^2_\gamma}^2 + \|d^* u\|_{L^2_\gamma}^2 
                &= \int_{\Omega} \Big[ |\nabla \mathbf{u}|^2 + (\operatorname{div}\mathbf{u})^2 - \sum_{i,j} \partial_i u_j \partial_j u_i \nonumber \\
                &\qquad - 2(\operatorname{div}\mathbf{u})(x\cdot \mathbf{u}) + (x\cdot \mathbf{u})^2 \Big]\,d\gamma. \label{eq:expanded}
            \end{align}
            We now calculate the two cross terms
            $$
            I_1 := -\int_{\Omega} \sum_{i,j} \partial_i u_j \partial_j u_i \,d\gamma, \qquad
            I_2 := -2\int_{\Omega} (\operatorname{div}\mathbf{u})(x\cdot \mathbf{u})\,d\gamma.
            $$
            Let $\partial_k^* = \partial_k - x_k$ for each $k=1,\ldots,n$. Then for each fixed pair $(i,j)$,
            \begin{align*}
                \int_{\Omega} \partial_i u_j \partial_j u_i \,d\gamma 
                &= -\int_{\Omega} u_i \partial_j^*(\partial_i u_j )\,d\gamma+\int_{\partial\Omega} u_i\partial_i u_j\nu_j\,d\gamma_\sigma  \\
                &= -\int_{\Omega} u_i \big( \partial_j\partial_i u_j \big)\,d\gamma + \int_{\Omega} u_i \partial_i u_j\,x_j\,d\gamma+\int_{\partial\Omega} u_i\partial_i u_j\nu_j\,d\gamma_\sigma.
            \end{align*}
            Summing over all $i,j$ yields
            \begin{align}
                I_1 &= \int_{\Omega} \mathbf{u} \cdot\nabla(\operatorname{div}\mathbf{u})\,d\gamma  -  \sum_{i,j}\int_{\Omega} u_i x_j\partial_i u_j\, d\gamma-\sum_{i,j}\int_{\partial\Omega}u_i\partial_i u_j\nu_j\,d\gamma_\sigma. \label{eq:I1}
            \end{align}
            Similarly, we have
            \begin{align}
                I_2 = 2\int_{\Omega} |\mathbf{u}|^2\,d\gamma + 2\sum_{i,j}\int_{\Omega}  u_i x_j \partial_i u_j\,d\gamma - 2\int_{\Omega} (x\cdot \mathbf{u})^2\,d\gamma - 2\int_{\partial\Omega}(x\cdot \mathbf{u})\as{\nu,u}\,d\gamma_\sigma. \label{eq:I2}
            \end{align}
            Adding \eqref{eq:I1} and \eqref{eq:I2}, we obtain
            \begin{align*}
                I_1 + I_2 
                &= \int_{\Omega} \mathbf{u} \cdot\nabla(\operatorname{div}\mathbf{u})\,d\gamma  
                + \sum_{i,j}\int_{\Omega}  u_i x_j \partial_i u_j\,d\gamma + 2\int_{\Omega} |\mathbf{u}|^2\,d\gamma\\
                &\quad  - 2\int_{\Omega} (x\cdot \mathbf{u})^2\,d\gamma -\sum_{i,j}\int_{\partial\Omega}u_i\partial_i u_j\nu_j\,d\gamma_\sigma 
                - 2\int_{\partial\Omega}(x\cdot \mathbf{u})\as{\nu,u}\,d\gamma_\sigma.
            \end{align*}
            Integration by parts gives
            \begin{align*}
                \int_{\Omega} \mathbf{u} \cdot\nabla(\operatorname{div}\mathbf{u})\,d\gamma 
                =  -\int_{\Omega} (\operatorname{div}\mathbf{u})^2\,d\gamma + \int_{\Omega} (\operatorname{div}\mathbf{u})(x\cdot \mathbf{u})\,d\gamma + \int_{\partial\Omega}(\operatorname{div}\mathbf{u})\as{\nu,u}\,d\gamma_\sigma
            \end{align*}
            and
            $$
            \sum_{i,j}\int_{\Omega}  u_i x_j \partial_i u_j\,d\gamma 
            = -\int_{\Omega} (\operatorname{div}\mathbf{u})(x\cdot \mathbf{u})\,d\gamma - \int_{\Omega} |\mathbf{u}|^2\,d\gamma + \int_{\Omega} (x\cdot \mathbf{u})^2\,d\gamma + \int_{\partial\Omega}(x\cdot \mathbf{u})\as{\nu,u}\,d\gamma_\sigma.
            $$
            Substituting these two equations into $I_1+I_2$ gives
            \begin{align*}
                I_1 + I_2 
                &= -\int_{\Omega} (\operatorname{div}\mathbf{u})^2\,d\gamma + \int_{\Omega} |\mathbf{u}|^2\,d\gamma - \int_{\Omega} (x\cdot \mathbf{u})^2\,d\gamma\\
                &\quad  -\sum_{i,j}\int_{\partial\Omega}u_i\partial_i u_j\nu_j\,d\gamma_\sigma 
                -\int_{\partial\Omega}(x\cdot \mathbf{u})\as{\nu,u}\,d\gamma_\sigma
                + \int_{\partial\Omega}(\operatorname{div}\mathbf{u})\as{\nu,u}\,d\gamma_\sigma.
            \end{align*}
                A direct calculation shows
            \begin{align*}
                \sum_{i,j=1}^n u_i \, \partial_i u_j \, \nu_j 
            = \mathbf{u}^\top \cdot \nabla^\top \as{\nu,u} 
            + \as{\nu,u} \, \partial_\mathbf{n} \as{\nu,u} 
            -\mathrm{II}(\mathbf{u}^\top, \mathbf{u}^\top),
            \end{align*}
            where $\nabla^\top f = \nabla f-(\partial_\mathbf{n} f)\mathbf{n}$ is the tangential differential. 

            A direct calculation shows that $\nu \vee u=\as{\nu, u}$. Therefore,
            $$
            u_{\mathrm{tan}} = u - \langle \nu, u \rangle \nu,\qquad
            u_{\mathrm{nor}} = \langle \nu, u \rangle \nu.
            $$
            Note that the condition $u_{\mathrm{tan}} = u - \langle \nu, u \rangle \nu=0$ means $\mathbf{u}^\top=0$ on $\Sigma$, and the condition  $u_{\mathrm{nor}} = \as{\nu, u} \nu=0$ means $\as{\nu, u}=0$ on $\Sigma^\prime$. Therefore, 
            \begin{align*}
                \sum_{i,j}\int_{\partial\Omega}u_i\partial_i u_j\nu_j\,d\gamma_\sigma 
                &= \frac{1}{2}\int_{\Sigma}\partial_\mathbf{n} \as{\nu,u}^2\,d\gamma_\sigma 
            - \int_{\Sigma^\prime}\mathrm{II}(\mathbf{u}^\top, \mathbf{u}^\top)\,d\gamma_\sigma. 
            \end{align*}
            Since $\Sigma\subset\overline{\Gamma}_2\cap\{x_k=0\}$ and $\nu = -dx_k$ on $\Sigma$ for some $k$, we have $x\cdot \mathbf{u}=0$ on $\Sigma$ and thus $(x\cdot \mathbf{u})\as{\nu,u}=0$ on $\partial\Omega$. Hence,
            \begin{align*}
                \int_{\partial\Omega}(x\cdot \mathbf{u})\as{\nu,u}\,d\gamma_\sigma=0.
            \end{align*}
            Similarly, for  $\Sigma\subset\overline{\Gamma}_2\cap\{x_k=0\}$, we have     $\operatorname{div}\mathbf{u} = \partial_k u_k = \partial_\mathbf{n}\as{\nu,u}$. Hence,
            \begin{align*}
                \int_{\partial\Omega}(\operatorname{div}\mathbf{u})\as{\nu,u}\,d\gamma_\sigma 
            = \frac{1}{2}\int_{\Sigma}\partial_\mathbf{n} \as{\nu,u}^2\,d\gamma_\sigma.
            \end{align*}
            Therefore, from equation \eqref{eq:expanded}, we have
            \begin{align}\label{E_second_fundamental}
                \|du\|^2_{L^2_\gamma} + \|d^* u\|^2_{L^2_\gamma}
                &= \int_{\Omega} \Big( |\nabla \mathbf{u}|^2 + (\operatorname{div}\mathbf{u})^2 \Big)\,d\gamma+ I_1+I_2  \notag\\
                &= \int_{\Omega} \big( |\nabla \mathbf{u}|^2 + |\mathbf{u}|^2 \big)\,d\gamma+ \int_{\Sigma^\prime}\mathrm{II}(\mathbf{u}^\top, \mathbf{u}^\top)\,d\gamma_\sigma\\
                &= \int_{\Omega} \big( |\nabla \mathbf{u}|^2 + |\mathbf{u}|^2 \big)\,d\gamma
                + \int_{\partial\Omega}\mathrm{II}(\mathbf{u}^\top, \mathbf{u}^\top)\,d\gamma_\sigma,\notag
            \end{align}
            which gives the equality \eqref{E:eta_u_equivalent} when $u\in C^\infty(\Omega,\Lambda^1)$. For $u\in H_\gamma(\Omega, \Sigma, d) \cap H_\gamma(\Omega,\Sigma^\prime,d^*)$, the equality \eqref{E:eta_u_equivalent} follows from the trace theorem and the density of the subspace $C^\infty(\Omega,\Lambda^1)\cap H_\gamma(\Omega, \Sigma, d) \cap H_\gamma(\Omega,\Sigma^\prime,d^*)$ in $H_\gamma(\Omega, \Sigma, d) \cap H_\gamma(\Omega,\Sigma^\prime,d^*)$.

            Now we want to show $\av{u}\in H_\gamma(\Omega, \Sigma, d) \cap H_\gamma(\Omega,\Sigma^\prime,d^*)$. That is, $\av{u}_{tan}=0$ on $\Sigma$ and $\av{u}_{nor}=0$ on $\Sigma^\prime$. 

        	On $\Sigma$, since ${u}_{tan}=0$, we have $\nu\wedge{u}=0$. By Cartan's Lemma, if we write $\nu=\nu_{k}dx_{k}$, then
        	\begin{align}\label{equationSigma}
        		u = c\nu = c\nu_{k}dx_{k}.
        	\end{align}
        	Thus $\nu\wedge\av{u}=0$. That is, $\av{u}_{tan}=0$.
        	
        	On $\Sigma^\prime$, we write $\nu=\nu_{j_1}dx_{j_1}+\nu_{j_2}dx_{j_2}$. Then the condition  $u_{\mathrm{nor}} = \as{\nu, u} \nu=0$ means 
        	\begin{align*}
        		\nu_{j_1}u_{j_1}+\nu_{j_2}u_{j_2}=0.
        	\end{align*} 
        	Since $\nu_{j_1}\cdot\nu_{j_2}\leq 0$ and $\nu_{j_1}^2+\nu_{j_2}^2=1$ on $\Sigma^\prime$, we know that $u_{j_1}\cdot u_{j_2}\geq 0$. Thus, $\av{u_{j_1}}\cdot \av{u_{j_2}}=u_{j_1}\cdot u_{j_2}$ and 
        	\begin{align*}
        		\nu\wedge*\av{u} = (\nu_{j_1}\av{u_{j_1}}+\nu_{j_2}\av{u_{j_2}})dx_1\wedge\cdots\wedge dx_n = 0,
        	\end{align*}
        	which implies $\av{u}_{nor}=0$ on $\Sigma^\prime$. Combining the above results, we have $\av{u}\in H_\gamma(\Omega, \Sigma, d) \cap H_\gamma(\Omega,\Sigma^\prime,d^*)$. 

            By Lemma \ref{regularpoints}, we may assume that for each fixed point $p=(p_1,\ldots,p_n)$ and\, for any $x\in \mathcal{O}(p)$, $\nu(x) = \nu_{n-1}dx_{n-1} + \nu_n dx_n$ with $\nu_n\neq 0$. Assume that in $\mathcal{O}(p)$ the boundary can be written as the graph of a smooth function $F:U\to\mathbb{R}$, where $U\subset\mathbb{R}^{n-1}$ is open. Precisely, we have the parametrization  
            $$
            \Phi(x_1,\dots,x_{n-1}) = (x_1,\dots,x_{n-1},\,F(x_1,\dots,x_{n-1})).
            $$
            The natural tangent vectors are  
            $$
            X_i := \partial_i\Phi = e_i + (\partial_i F)\,e_n, \qquad i=1,\dots,n-1,
            $$
            where $e_1,\dots,e_n$ is the standard basis of $\mathbb{R}^n$. Then we obtain
            \begin{align*}
                \langle \nu, X_i\rangle = \nu_{n-1}\langle e_{n-1}, X_i\rangle +\nu_{n}\langle e_n, X_i\rangle = 0,
            \end{align*}
            which implies $\partial_i F = 0$  for each $i=1,\dots,n-2$ and $\partial_{n-1} F = -\nu_{n-1}/\nu_n$. Therefore, $F$ depends only on $x_{n-1}$ and $\nu$ depends only on $x_{n-1}$ and $x_{n}$. The boundary condition  $\as{\nu, u}=0$ on $\Sigma^\prime$ forces $u_{n-1} \nu_{n-1} + u_n \nu_{n} = 0$, so we may write $u_{n-1} = t\nu_{n}, u_n = -t\nu_{n-1}$,where $t$ is a scalar satisfying $t^2 = |\mathbf{u}|^2$. Using the identity $\nu_{n-1} \partial_i\nu_{n-1} + \nu_{n} \partial_i\nu_{n} = 0$, we can derive
            $$
            \partial_i\nu_{n-1} = a_i \nu_{n}, \qquad \partial_i\nu_{n} = -a_i\nu_{n-1}
            $$
            for $i=1,\ldots,n$ and some scalar $a_i$. Hence,
            \begin{align*}
                \mathrm{II}(\mathbf{u}^\top, \mathbf{u}^\top)
                =\sum_{i,j} u_i (\partial_i \nu_j) u_j  = t^2(a_{n-1} \nu_{n} - a_n \nu_{n-1})    = \kappa\, |\mathbf{u}|^2 ,
            \end{align*}
            where $\kappa = a_{n-1} \nu_{n} - a_n \nu_{n-1}$. Direct computation shows that $\kappa$ is precisely the unique possible non‑zero principal curvature at the point. Moreover, $\kappa$ is nonnegative whenever $\Omega$ is convex. 

            Therefore, it follows from the equality \eqref{E_second_fundamental} that
            \begin{align*}
                \|d\av{u}\|^2_{L^2_\gamma} + \|d^* \av{u}\|^2_{L^2_\gamma}\leq \|du\|^2_{L^2_\gamma} + \|d^* u\|^2_{L^2_\gamma},
            \end{align*}
            which completes the proof.
        \end{proof}
        As a consequence of the above lemma, we have the following key lemma.
        \begin{lemma}\label{abs}
        	Let $ \Omega\subset \mathcal{M}^n$ be a lip domain satisfying Assumption \ref{assumption}, and let $\Sigma=\mathrm{int}\LRc{{\Gamma}_2\cap\{x_k=0\}}$ be a connected and open subset for some integer $1\leq k\leq n$. Assume that $\Omega$ is creased relative to $\Sigma$. Then for any minimizer $u\in H_\gamma(\Omega, \Sigma, d) \cap H_\gamma(\Omega,\Sigma^\prime,d^*)$ of the functional (\ref{variation}), $\av{u}$ is also a minimizer.
        \end{lemma}
        
        The following lemma is an extension of results in \cite{KennedyJamesB2024Oths} to Gaussian settings.
        \begin{lemma}\label{u>0}
        	Let $ \Omega\subset \mathcal{M}^n$ be a lip domain satisfying Assumption \ref{assumption}, and let $\Sigma=\mathrm{int}\LRc{{\Gamma}_2\cap\{x_k=0\}}$ be a connected and open subset for some integer $1\leq k\leq n$. Assume that $\Omega$ is creased relative to $\Sigma$. Assume that a solution $u$ of the equation (\ref{$1$-form}) satisfies $u_j\geq 0$ in $\Omega$ or $u_j\leq 0$ in $\Omega$ for all $j=1,\ldots,n$. If $d^* u=0$, then $u\equiv 0$.
        \end{lemma}
        \begin{proof}
        	If $\Sigma=\emptyset$, then for each $j\in\{1,\ldots,n\}$ such that $u_j\geq 0$, by Stokes' theorem, 
        	\begin{align}\label{u_equiv_0}
        		0\leq\int_{\Omega}\as{dx_j, u}d\gamma = \int_{\Omega}\as{x_j, d^* u}d\gamma + \int_{\partial\Omega}x_j\as{\nu, u}d\gamma_{\sigma}=0,
        	\end{align}
        	which implies $u_j\equiv 0$. Similarly, we also have $u_j\equiv 0$ if $u_j\leq 0$ in $\Omega$. Thus, $u\equiv 0$ in $\Omega$. 
        	
        	Now we assume $\Sigma\neq \emptyset$. 
        	Since $\Sigma=\mathrm{int}\LRc{{\Gamma}_2\cap\{x_k=0\}}$, we may assume that for each point $p\in\Sigma$, 
        	\begin{align*}
        		\nu=-dx_k.
        	\end{align*}
        	
        	For $f(x):=\mathrm{dist}(x,\Sigma)=x_k$, we have $f|_\Sigma=0, f>0$ and  $df(x)=dx_k$ in $\Omega$. We may assume $u_k\geq 0$ in $\Omega$. Thus, by Stokes' theorem, we have
        	\begin{align*}
        		0\leq\int_{\Omega}\as{df, u}d\gamma = \int_{\Omega}\as{f, d^* u}d\gamma + \int_{\Sigma}f\as{\nu, u}d\gamma_{\sigma}+\int_{\Sigma^\prime}f\as{\nu, u}d\gamma_{\sigma}=0,
        	\end{align*}
        	which implies 
        	\begin{align*}
        		\as{df, u}=\as{dx_k, u}\equiv 0\quad\text{in}\ \Omega.
        	\end{align*}
        	Thus, $\as{\nu, u}\equiv 0$ on $\partial\Omega$. It follows from the inequality (\ref{u_equiv_0}) that $u_j\equiv 0$ in $\Omega$ for each $j=1,\ldots,n$. That is, $u\equiv 0$ in $\Omega$.
        \end{proof}

        Using the analyticity of the eigenforms of problem \eqref{$1$-form}, we now prove that for every first nontrivial eigenfunction $\psi$ with mixed boundary condition, $d\psi$ is a minimizer of the functional (\ref{variation}). Moreover, any minimizer of the functional (\ref{variation}) is of the form $d\psi$ for some first nontrivial eigenfunction $\psi$.
        \begin{lemma}\label{mu-2}
        	Let $ \Omega\subset \mathcal{M}^n$ be a lip domain satisfying Assumption \ref{assumption}, and let $\Sigma=\mathrm{int}\LRc{{\Gamma}_2\cap\{x_k=0\}}$ be a connected and open subset for some integer $1\leq k\leq n$. Assume that $\Omega$ is creased relative to $\Sigma$. If $\Sigma=\emptyset$, then $\eta_1=\mu_2$. If $\Sigma\neq\emptyset$, then $\eta_1=\lambda_1^{\Sigma}$. In particular, if $u$ is a minimizer of the functional (\ref{variation}), then $u\in dH_\gamma(\Omega, \Sigma, d)$ and the following trichotomy holds.
        	\begin{align}\label{trichotomy}
        		u_j\geq 0\ \text{in}\ \Omega,\ u_j\leq 0 \ \text{in}\  \Omega,\ \text{or}\ u_j\equiv 0 \ \text{in}\  \Omega,\qquad j=1,\ldots,n.
        	\end{align}
        \end{lemma}
        \begin{proof}
        	Consider the case $\Sigma\neq\emptyset$ (the proof for $\Sigma=\emptyset$ is similar). Assume by contradiction that $\eta_1 <\lambda_1^{\Sigma}$. Then by Hodge decomposition (\ref{Hodgedecomposition}), there is some minimizer $u\in H_\gamma(\Omega, \Sigma, d) \cap H_\gamma(\Omega,\Sigma^\prime,d^*)$ of the functional (\ref{variation}) such that $d^* u=0$. By Lemma \ref{abs}, we know that $\av{u}$ is also a minimizer.  Since $\eta_1 <\lambda_1^{\Sigma}$, we have $\av{u}\geq 0$ and $d^* \av{u}=0$. Thus, $u$ satisfies the conditions of Lemma \ref{u>0}. Therefore, $\av{u}\equiv 0$ and $u\equiv 0$, which is a contradiction. Hence, we obtain 
        	\begin{align*}
        		\eta_1=\lambda_1^{\Sigma}.
        	\end{align*}
        	Note that $\av{u}$ is a solution of the equation (\ref{$1$-form}) with $\eta=\eta_1$. From the regularity theory, we know that $u$ and $\av{u}$ are analytic. Thus, if $u_i(x)>0$ for some $i\in\{1,\ldots,n\}$ and some point $x\in\Omega$, then there is an open neighborhood $U_0$ of $x$ in which $u_i>0$. So we get $u_i=\av{u_i}$ in $U_0$. By the regularity of $u$ and $\av{u}$, we get
        	\begin{align*}
        		u_i\equiv\av{u_i}\geq 0 \quad\text{in}\ \Omega.
        	\end{align*}
        	Similarly, if $u_l(y)<0$ for some $l\in\{1,\ldots,n\}$ and some point $y\in\Omega$, then 
        	\begin{align*}
        		u_l\equiv -\av{u_l}\leq 0 \quad\text{in}\ \Omega.
        	\end{align*}
        	Therefore, for each $j\in\{1,\ldots,n\}$, we have the trichotomy \eqref{trichotomy}
. 

            If $u\notin dH_\gamma(\Omega, \Sigma, d)$ is a minimizer, by Hodge decomposition (\ref{Hodgedecomposition}),
        	\begin{align*}
        		u\in dH_\gamma(\Omega, \Sigma, d) \oplus d^* H_\gamma(\Omega, \Sigma', d^*) \oplus \mathscr{H}_\gamma(\Omega, \Sigma, \Sigma'),
        	\end{align*}
        	we may assume $d^* u = 0$ by subtracting terms in $dH_\gamma(\Omega, \Sigma, d)$. Then by Lemma \ref{u>0} we know that $u\equiv 0$ in $\Omega$, which leads to a contradiction. 
        	
        	Therefore, any minimizer of the functional (\ref{variation}) must be in $ dH_\gamma(\Omega, \Sigma, d)$. This proves the lemma.
        \end{proof}
        \begin{remark}
        	When $n=2$, as a corollary of Lemma \ref{mu-2}, for lip domains with $\Sigma=\emptyset$, we obtain the eigenvalue inequality that 
        	\begin{align*}
        		\mu_2< \lambda_1
        	\end{align*}
        	for the first nontrivial Neumann eigenvalue $\mu_2$ and the first Dirichlet eigenvalue $\lambda_1$;
        	and when $\Sigma\neq\emptyset$, we have  
        	\begin{align}\label{eigeninequality}
        		\lambda_1^{\Sigma}< \lambda_1^{\Sigma^\prime}
        	\end{align}
        	for the first mixed eigenvalue with Dirichlet boundary on $\Sigma$ and the first mixed eigenvalue with Dirichlet boundary on $\Sigma^\prime$, where $\Sigma$ satisfies the conditions in Lemma \ref{u>0}. It is noteworthy that $\lambda_1^{\Sigma}< \lambda_1^{\Sigma^\prime}$ is a generalization of the results in \cite{MR4860866}, which deals with the case in $\mathbb{R}^2$. 
        	
        	In contrast, the strict inequality (\ref{eigeninequality}) fails when $\Sigma$ is not regular. For instance, we may take
        	\begin{align*}
        		\Omega=(-\pi,\pi)\times(-\pi,\pi),\ \Sigma=(\{-\pi\}\times (-\pi,\pi))\cup((-\pi,\pi)\times\{-\pi\}),
        	\end{align*} 
        	then $\lambda_1^{\Sigma}=\lambda_1^{\Sigma^\prime}$ for $\Sigma^\prime=\partial\Omega\backslash\overline{\Sigma}$. 

            Note that if $\Sigma$ is not connected, then
        \begin{align*}
        	dH_{\Sigma,\gamma}^0:=\{d u\in L^2_\gamma(\Omega,\Lambda^1):u\in L^2_\gamma(\Omega), u=0\ \text{on}\ \Sigma\}\varsubsetneqq H^1_\gamma(\Omega, \Sigma, d)\subset L^2_\gamma(\Omega,\Lambda^1),
        \end{align*} 
        which implies that the minimizer of (\ref{variation}) may not be contained in $dH_{\Sigma,\gamma}^0$ and that the first mixed eigenfunction may fail to be monotonic in certain coordinate directions $e_i$. Thus, the connectivity of $\Sigma$ is important in our theorems.
        \end{remark}
        
        Now we are ready to prove the hot spots conjecture for lip domains.
        
        ~\\
        \emph{Proof of Theorem \ref{Tlip}}.
        Assume $\psi$ is the first nontrivial eigenfunction of $\Omega$ with mixed boundary conditions, associated with the first nontrivial mixed eigenvalue $\lambda_1^{\Sigma}$. By Lemma \ref{mu-2}, $u=d\psi$ is a minimizer of the functional (\ref{variation}) and the trichotomy (\ref{trichotomy}) holds. That is, $u_j$ does not change its sign in $\Omega$ for all $j=1,\ldots,n$. Then we have the trichotomy for each $j=1,\ldots,n$,
        \begin{align*}
        	\partial_j\psi\geq 0\ \text{in}\ \Omega,\ 	\partial_j\psi\leq 0 \ \text{in}\  \Omega,\ \text{or}\ 	\partial_j\psi\equiv 0 \ \text{in}\  \Omega,
        \end{align*}
        which implies $\psi$ is monotonic in each $e_j$ direction. If
        \begin{align*}
        	\psi(p)=\max_{\Omega}\psi
        \end{align*}
        for some interior point $p\in\Omega$, then for each point $q\in\Omega$ such that $q_j-p_j>0$ where $	\partial_j\psi\geq 0\ \text{in}\ \Omega$ and $q_j-p_j<0$ where $	\partial_j\psi\leq 0\ \text{in}\ \Omega$, we have 
        \begin{align*}
        	\psi(q)=\psi(p).
        \end{align*} 
        The set of such $q$ contains a nonempty open subset of $\Omega$, and by the unique continuation principle for $\psi$, we have $\psi\equiv\psi(p)$ in $\Omega$, which leads to a contradiction. Thus, the extrema of $\psi$ are attained only on the boundary.
        \hfill$\square$~\\\par
        \begin{remark}
        	If $\Sigma=\emptyset$, then $\lambda_1^{\Sigma}=\mu_2$. Therefore, every first nontrivial Neumann eigenfunction, as well as every first mixed eigenfunction, is monotonic in the lip domain in each $e_i$ direction for $i=1,\ldots, n$.
        \end{remark}

        \section{\texorpdfstring{$n$}{n}-Symmetric Domains Have No Hot Spots}\label{sec5}
        In this section, we begin with the following lemma to reduce the case of symmetric domains to lip domains.
        
        \begin{lemma}\label{n-symmetric}
        	Let $\Omega\subset \mathcal{M}^n$ be a bounded  $n$-symmetric domain satisfying Assumption \ref{assumption}.  Then for any $k\in\{1,\ldots,n\}$, any first $k$-antisymmetric Neumann eigenfunction $ w_k $ is unique and symmetric with respect to the other coordinate hyperplane $ \{x_{k^\prime} = 0\},\ \forall\ k^\prime \neq k $.
        \end{lemma}
        \begin{proof}
        	Assume $k=n$. First, note that every Neumann eigenfunction $w$ in $\Omega$, which is antisymmetric with respect to reflection in the hyperplane $ \{x_n = 0\} $, corresponds to an eigenfunction $w$ in $\Omega\cap\{x_n > 0\}$ with Dirichlet boundary condition on $\Omega\cap\{x_n = 0\}$ and Neumann boundary condition on $\partial\Omega\cap\{x_n > 0\}$. Thus, we only need to consider $w$ in $\Omega\cap\{x_n > 0\}$.  From the condition that $w$ is the first $n$-antisymmetric Neumann eigenfunction, we know that $w$ is the first eigenfunction with mixed boundary condition on $\Omega\cap\{x_n > 0\}$. Hence, $w$ is unique and does not change sign. After possibly replacing $w$ by $-w$,  we may assume $w>0$ in $\Omega\cap\{x_n > 0\}$. 
        	
        	Let 
        	\begin{align*}
        		\widetilde{w}(x_1,x_2,\ldots, x_n):=w(-x_1,x_2,\ldots, x_n). 
        	\end{align*}
        	Then $\widetilde{w}$ is also a first eigenfunction with mixed boundary condition in $\Omega\cap\{x_n > 0\}$. From the uniqueness of $w$, we know that $w=\widetilde{w}$. Similarly, for any $k^\prime \neq n $, $w$ is symmetric with respect to the coordinate hyperplane $ \{x_{k^\prime} = 0\}$.
        \end{proof}
        
        Now we prove our main theorem.
        
        ~\\\emph{Proof of Theorem \ref{Tsymmetric}}. 
        Assume $k=n$. By Lemma \ref{n-symmetric}, we only need to consider $w_n$ in the lip domain
        \begin{align*}
        	\widetilde{\Omega}:=\Omega\cap\mathcal{O}
        \end{align*}
        with 
        \begin{align*}
        	w_n\Big|_{\partial\widetilde{\Omega}\cap\{x_n = 0\}}=0,\ \partiald{w_n}{\mathbf{n}}\Big|_{\partial\widetilde{\Omega}\cap\{x_n > 0\}}=0.
        \end{align*}
        Then this is a special case of Theorem \ref{Tlip} when $\Sigma=\mathrm{int}\LRc{\partial\widetilde{\Omega}\cap\{x_n = 0\}}$. Thus, $w_n$ is monotonic in $\widetilde{\Omega}$ in each $e_i$ direction for $i=1,\ldots,n$ and $w_n$ cannot have any interior maximum point except for points on some hyperplane. Since $\partial_n w_n>0$ on ${\partial\widetilde{\Omega}\cap\{x_n = 0\}}$ by Hopf's Lemma (after multiplication by a constant), we have $\partial_n w_n\geq 0$ in $\widetilde{\Omega}$.
        
        Assume by contradiction that for some $k\in\{1,\ldots,n-1\}$ and some $p\in\Omega\cap\{x_k = 0\}$,
        \begin{align*}
        	w_n(p)=\max_{\Omega}w_n.
        \end{align*}
        
        By symmetry and continuity of $\partial_n w_n$, we get $\partial_n w_n\geq 0$ in ${\Omega}$. Thus, we have
         \begin{align*}
         	w_n(p+te_n)=w_n(p)>0,\ \forall t>0, p+te_n\in\Omega.
         \end{align*}
         
         Since $w_n$ is analytic in $\Omega$, we have 
         \begin{align*}
         	w_n(p+te_n)\equiv w_n(p)>0,\ \forall t\in\mathbb{R}, p+te_n\in\Omega.
         \end{align*}
         But this contradicts the condition that $w_n=0$ on $\Omega\cap\{x_n = 0\}$. Thus, $w_n$ cannot have any interior maximum point in $\Omega$.      
        \hfill$\square$~\\\par
        \begin{remark}
            When $n=2$, Theorem \ref{Tsymmetric} provides a Gaussian analogue of  \cite{MR1775736}*{Theorem 1.1}. Note that if $\Omega$ is convex, then $\Omega\cap\{x_1<0,x_2>0\}$ is a lip domain and creased relative to $\partial\Omega\cap\{x_1<0,x_2=0\}$ and $\partial\Omega\cap\{x_1=0,x_2>0\}$. Therefore, if $\Omega$ is convex and the second Neumann eigenvalue $\mu_2$ is simple, then Theorem \ref{Tsymmetric} immediately implies Corollary \ref{Corollary} due to the symmetry of the domain. However, when $\mu_2$ is not simple, then there exist two linearly independent Neumann eigenfunctions. Even though these eigenfunctions can be taken to be $k$-antisymmetric, their linear combinations are no longer symmetric or antisymmetric with respect to the axes. Therefore, Theorem \ref{Tsymmetric} cannot be applied to derive Theorem 1.4 in \cite{MR1775736}, which asserts that for any first nontrivial eigenfunction $u$, there exists a direction $e \in \mathbb{R}^2 \setminus \{0\}$ such that $e \cdot \nabla u > 0$ in $\Omega$.
        \end{remark}

        Unfortunately, the Hodge decomposition and the variational method do not extend to lip domains in constant curvature spaces directly, because the inequality \eqref{E:absolute_minimizer} no longer holds there. Indeed, if Theorem \ref{Tlip} were true for lip domains in constant curvature spaces, then by symmetry of rectangles one could conclude that some first nontrivial Neumann eigenfunction depends on only one variable, which would contradict the proposition proved below. Hence, the behavior of the partial derivatives of eigenfunctions in constant curvature spaces genuinely differs from that in Euclidean spaces.
        
        \begin{proposition}\label{prop}
        	Let $\mathbb{D}$ be the Poincar\'{e} disk model of the hyperbolic plane and let $\Omega\subset\mathbb{D}$ be a rectangle centered at $0$. There is no first nontrivial Neumann eigenfunction $\psi$ on $\Omega$  such that $\psi$ depends only on $x_i$ for some $i\in\{1,2\}$.
        \end{proposition}
        \begin{proof}
        	The Hodge Laplacian $\Delta := d\delta + \delta d$ on $0$-forms is given by
        	\begin{align*}
        		\Delta \psi = -\frac{\LRc{1-\av{x}^2}^2\LRc{\partial_{11}\psi+\partial_{22}\psi}}{4}.
        	\end{align*}
            
        	We may assume $\partial_1\psi = 0$. Then
        	\begin{align*}
        		\Delta \psi =  -\frac{\LRc{1-\av{x}^2}^2\partial_{22}\psi}{4}=\mu_2\psi.
        	\end{align*}
        	Differentiating both sides with respect to $x_1$, we have
        	\begin{align*}
        		{\LRc{1-\av{x}^2}x_1\partial_{22}\psi}=0.
        	\end{align*}
        	Thus, $	\Delta \psi \equiv 0$ and $\mu_2\psi\equiv 0$ in $\Omega$, which means $\psi\equiv 0$ in $\Omega$. This is a contradiction.
        \end{proof}

	\section*{Acknowledgments}
    We would like to thank Prof. Zhiqin Lu for the discussions on differential forms, and also express our gratitude to Florentin M\"{u}nch and Haohang Zhang for their insightful discussions and suggestions.

	\bibliographystyle{plain}
	\bibliography{Hot_spot_conjecture_in_Gaussian_spaces}                        
	
\end{document}